\documentclass[11pt,a4paper]{article}
\usepackage{a4wide,amsfonts,amsmath,latexsym,amssymb,euscript,eufrak,graphicx,units,mathrsfs}

\usepackage{stmaryrd}
\usepackage{float}
\newfloat{figure}{H}{lof}
\floatname{figure}{\figurename}

\DeclareMathAlphabet{\eufrak}{U}{}{}{}  % Euler fraktur math
\SetMathAlphabet\eufrak{normal}{U}{euf}{m}{n}
\SetMathAlphabet\eufrak{bold}{U}{euf}{b}{n}

\usepackage{amsmath, amsthm, amsfonts, amssymb}

\numberwithin{equation}{section}

\def\real{{\mathord{{\rm I\kern-2.8pt R}}}}        % Fake blackboard bold R.
\def\inte{{\mathord{{\rm I\kern-2.8pt N}}}}
\def\PP{{\mathord{{\rm I\kern-2.8pt P}}}}

\def\real{{\mathord{\mathbb R}}}

\def\inte{{\mathord{\mathbb N}}}

\def\Var{{\mathrm{{\rm Var}}}}

\def\HH{\EuFrak H}
\def\NN{\mathscr N}

\newcommand{\ind}{\mathbf{1}}

\def\ee{E}

\def\E{\mathop{\hbox{\rm I\kern-0.20em E}}\nolimits}
\def\Var{\mathop{\hbox{\rm Var}}\nolimits}

\def\real{\mathbb{R}}

\newtheorem{prop}{Proposition}[section]

\newtheorem{theorem}[prop]{Theorem}

\allowdisplaybreaks

\begin{document}

\begin{center}
{\large \textbf{Error bounds on the non-normal approximation of Hermite power variations of fractional Brownian motion}}\\[0pt]
~\\[0pt]
Jean-Christophe Breton%
\footnote{%
Universit\'e de La Rochelle, Laboratoire Math\'ematiques, Image et Applications, Avenue Michel Cr\'epeau, 17042 La Rochelle Cedex, France, \texttt{jean-christophe.breton@univ-lr.fr}}
and Ivan Nourdin\footnote{%
Laboratoire de Probabilit\'es et Mod\`eles Al\'eatoires,
Universit\'e Pierre et Marie Curie (Paris VI), Bo\^ite courrier
188, 4 place Jussieu, 75252
Paris Cedex 05, France, \texttt{ivan.nourdin@upmc.fr}}
\\[0pt]
{\it Universit\'e de La Rochelle and Universit\'e Paris VI}\\
~\\[0pt]
\end{center}

{\small \noindent \textbf{Abstract:} 
Let $q\geq 2$ be a positive integer, $B$ be a fractional Brownian motion with Hurst index $H\in(0,1)$, 
$Z$ be an Hermite random variable of index $q$,
and $H_q$ denote the Hermite polynomial having degree $q$.
For any $n\geq 1$, set $V_n=\sum_{k=0}^{n-1} H_q(B_{k+1}-B_k)$.
The aim of the current paper is to derive, in the case when the Hurst index verifies $H>1-1/(2q)$, an upper bound for the total variation distance between the laws 
$\mathscr{L}(Z_n)$ and $\mathscr{L}(Z)$, where $Z_n$ stands for the correct renormalization of $V_n$
which converges in distribution towards $Z$. 
Our results
should be compared with those obtained recently by Nourdin and Peccati (2007) in the case when $H<1-1/(2q)$,
corresponding to the situation where one has normal approximation.\\
\normalsize
%\newline
}

{\small \noindent \textit{Key words:} Total variation distance;
Non-central limit theorem; 
Fractional Brownian motion; 
Hermite power variation; 
Multiple stochastic integrals;
Hermite random variable.~\newline
\normalsize }

{\small \noindent \textit{Current version}: April 2008}%\newline
\normalsize

\section{Introduction}
Let $q\geq 2$ be a positive integer and $B$ be a fractional Brownian motion (fBm) with Hurst index $H\in(0,1)$. 
The asymptotic behavior of the $q$-Hermite power variations of $B$ with respect to $\mathbb{N}$,  defined as
\begin{equation}\label{her-intro}
V_n=\sum_{k=0}^{n-1} H_q(B_{k+1}-B_k),\quad n\geq 1,
\end{equation}
has recently received a lot of attention,
see e.g. \cite{N}, \cite{NNT} and references therein. 
Here, $H_q$ stands for the Hermite polynomial with degree $q$, given by 
$
H_q(x)=(-1)^q e^{x^2/2}\frac{d^q}{dx^q}\big(e^{-x^2/2}\big).
$
We have $H_2(x)=x^2-1$, $H_3(x)=x^3-3x$, and so on. 
The analysis of the asymptotic behavior of (\ref{her-intro}) is motivated, for instance, by 
the traditional applications of quadratic variations to parameter estimation problems
(see e.g. \cite{begyn,coeurjolly,istas-lang,tudor-viens} and references therein).

In the particular case of the standard Brownian motion (that is when $H=\frac 12$), 
the asymptotic behavior of \eqref{her-intro}
can be immediately deduced from the classical central limit theorem. 
When $H\neq \frac12$,  the increments of $B$ are not  independent anymore and
the asymptotic behavior of \eqref{her-intro} is consequently more difficult to reach. 
However, thanks to the seminal works of
Breuer and Major \cite{BM}, Dobrushin and Major \cite{DoMa}, Giraitis and Surgailis \cite{GS} 
and Taqqu \cite{T}, it is well-known that we have,
as $n\to\infty$:
\begin{enumerate}
\item If $0<H<1-1/(2q)$ then
\begin{equation}\label{eq:Breuer_Major1}
Z_n:=\frac{V_n}{ \sigma_{q,H}\,\sqrt n}
\,\overset{{\rm Law}}{\longrightarrow}\,
\NN(0,1).
\end{equation}
\item If $H=1-1/(2q)$ then
\begin{equation}
\label{eq:Breuer_Major2}
Z_n:=\frac{V_n}{\sigma_{q,H}\sqrt{ n \log n}}
\,\overset{{\rm Law}}{\longrightarrow}\,
\NN(0,1).
\end{equation}
\item If $H>1-1/(2q)$ then
\begin{equation}
\label{eq:Breuer_Major3}
Z_n:=\frac{V_n}{n^{1-q(1-H)}} 
\,\overset{{\rm Law}}{\longrightarrow}\,
Z \sim\mbox{``Hermite random variable''.}
\end{equation}
\end{enumerate}
Here, $\sigma_{q,H}>0$ denotes an (explicit) constant depending only on $q$ and $H$. Moreover,
the Hermite random variable $Z$ appearing in (\ref{eq:Breuer_Major3}) is defined as the value at time
$1$ of the Hermite process, i.e. 
\begin{equation}\label{her-rv}
Z=I_q^W(L_1),
\end{equation} 
where $I_q^W$ denotes the $q$-multiple stochastic integral with respect to a Wiener process $W$,
while $L_1$ is the symmetric kernel defined as
$$
L_1(y_1, \dots, y_q)=\frac 1{q!}\ind_{[0,1]^q}(y_1, \dots, y_q) \int_{y_1\vee \dots \vee y_q}^1 
\partial_1 K_H(u,y_1)\dots \partial_1 K_H(u,y_q) du,
$$
with $K_H$ the square integrable kernel given by (\ref{kh}).
We refer to \cite{NNT} for a complete discussion of this subject.

The exact expression of the distribution function of $Z_n$ is very complicated when $H\neq \frac12$.
That is why, when $n$ is assumed to be large, it is common to use
\eqref{eq:Breuer_Major1}--\eqref{eq:Breuer_Major3} as a justification to replace, in any computation involving it, 
the distribution function of $Z_n$ by that of the corresponding limit.
Of course, if one applies this strategy without any care (in particular, if one has no idea, even imprecise, of an 
error bounds in terms of $n$), then it is easy to imagine that the obtained result could be very far from the reality
(as a ``concrete'' example, see (\ref{eq:H=}) below, which represents the worst case we will obtain here).
To the best of our knowledge, in all the works using \eqref{eq:Breuer_Major1}--\eqref{eq:Breuer_Major3}
with statistical applications in mind (for instance \cite{begyn,coeurjolly,istas-lang,tudor-viens}), 
never their author(s) considered this problem. The current paper, together with \cite{NP},
seem to be the first attempt in such direction. 

Recall that the {\sl total variation distance} between the laws of two real-valued
random variables $Y$ and $X$ is defined as
$$
d_{TV}\big(\mathscr{L}(Y),\mathscr{L}(X)\big)  = \sup_{A\in\mathscr{B}(\mathbb{R})} \big| P(Y\in A)-P(X\in A)\big|
$$
where $\mathscr{B}(\mathbb{R})$ denotes the set of Borelian of $\mathbb{R}$.
In \cite{NP}, by combining Stein's method with Malliavin calculus (see also Theorem \ref{theo:NP} below),
the following result is shown:
\begin{theorem}
\label{theo:NP_Berry}
If $H<1-1/(2q)$ then, for some constant $c_{q,H}>0$ depending uniquely on $q$ and $H$, we have:  
$$
d_{TV}\big(\mathscr{L}(Z_n),\NN(0,1)\big)
 \leq c_{q,H}
\left\{
\begin{array}{ll}
n^{-1/2} &\mbox{if } H\in (0, \frac12]\\
&\\
n^{H-1} &\mbox{if } H\in [\frac12,\frac{2q-3}{2q-2}]\\
&\\
n^{qH-q+\frac 12} &\mbox{if } H\in [\frac{2q-3}{2q-2}, 1- \frac1{2q})
\end{array}
\right.
$$
for $Z_n$ defined by (\ref{eq:Breuer_Major1}).
\end{theorem}
Here, we deal with the remaining cases, that is when $H\in [1-\frac{1}{2q}, 1)$. 
Our main result is as follows: 
\begin{theorem}
\label{theo:BN}
\begin{enumerate}
\item If $H=1-1/(2q)$ then, for some constant $c_{q,H}>0$ depending uniquely on $q$ and $H$, we have\begin{equation}
\label{eq:H=}
d_{TV}\big(\mathscr{L}(Z_n),\NN(0,1)\big)
\leq \frac{c_{q,H}}{\sqrt{\log n}}
\end{equation}
for $Z_n$ defined by (\ref{eq:Breuer_Major2}).

\item If $H\in (1-1/(2q), 1)$ then, for some constant $c_{q,H}>0$ depending uniquely on $q$ and $H$, we have
\begin{equation}
\label{eq:H>}
d_{TV}\big(\mathscr{L}(Z_n),\mathscr{L}(Z)\big)\leq c_{q,H}\, n^{1-\frac1{2q}-H}
\end{equation}
for $Z_n$ and $Z$ defined by (\ref{eq:Breuer_Major3}).
\end{enumerate}
\end{theorem}

Actually, the case when $H=1-1/(2q)$ can be tackled by mimicking the proof of Theorem \ref{theo:NP_Berry}. 
Only minor changes are required: we will also conclude thanks to the following general 
result by Nourdin and Peccati.
\begin{theorem}
\label{theo:NP} (cf. \cite{NP})
Fix an integer $q\geq 2$ and let 
$\{f_n\}_{n\geq 1}$ be a sequence of $\HH^{\odot q}$.
Then we have
$$
d_{TV}\left(\mathscr{L}\big(I_q(f_n)\big),\mathscr{N}(0,1)\right)
\leq 2\,\sqrt{E\left(
1-\frac1q\|DI_q(f_n)\|^2_{\HH}
\right)^2}, 
$$
where $D$ stands for the Malliavin derivative with respect to $X$.
\end{theorem}
Here, and for the rest of the paper, $X$ denotes a centered Gaussian isonormal process on 
a real separable Hilbert space $\HH$ and, as usual, $\HH^{\odot q}$ (resp. $I_q$) stands for the $q$th symmetric tensor product of $\HH$
(resp. the multiple Wiener-It\^o integral of order $q$ with respect to $X$). See Section 
\ref{sec:preliminaire} for more precise definitions and properties.

When $H\in (1-1/(2q),1)$, 
Theorem \ref{theo:NP} can not be used (the limit in (\ref{eq:Breuer_Major3}) being not Gaussian), 
and another argument is required. Our new idea is as follows. First, using the scaling property
(\ref{scaling}) of fBm, we construct, for every {\sl fixed} $n$, a copy $S_n$
of $Z_n$ that converges in $L^2$. Then, we use the following result by 
Davydov and Martynova.
\begin{theorem}
\label{theo:Davydov} (cf. \cite{DM}; see also \cite{B})
Fix an integer $q\geq 2$ and let $f\in \HH^{\odot q}\setminus\{0\}$. Then, 
for any sequence $\{f_n\}_{n\geq 1}\subset \HH^{\odot q}$ converging to $f$, their exists
a constant $c_{q,f}$, depending only on $q$ and $f$, such that:
$$
d_{TV}\left(\mathscr{L}\big(I_q(f_n)\big),\mathscr{L}\big(I_q(f)\big)\right) 
\leq c_{q,f} \|f_n-f\|_{\HH^{\odot q}}^{1/q}. 
$$
\end{theorem}

The rest of the paper is organized as follows. 
In Section \ref{sec:preliminaire}, some preliminary results on fractional Brownian motion and 
Malliavin calculus are presented. 
Section \ref{sec:H>} deals with the case $H\in (1-1/(2q), 1)$, 
while the critical case $H=1-1/(2q)$ is considered in Section \ref{sec:H=}. 

\section{Preliminaries}
\label{sec:preliminaire}

Let $B=\{B_t,\,t\geq 0\}$ be a fBm with Hurst index $H\in(0,1)$, 
that is a centered Gaussian process, started from zero and with covariance function 
$E(B_{s}B_{t})=R(s,t)$, where
$$
R(s,t)=\frac{1}{2}\left( t^{2H}+s^{2H}-|t-s|^{2H}\right);\quad
s,t\geq 0.
$$
In particular, it is immediately shown that $B$ has stationary increments and is selfsimilar of index $H$. 
Precisely,
for any $h,c>0$, we have
\begin{equation}\label{scaling}
\{B_{t+h}-B_h,\,t\geq 0\} \overset{{\rm Law}}{=} \{B_t,\,t\geq 0\}\quad\mbox{and}\quad
\{c^{-H}\,B_{ct},\,t\geq 0\} \overset{{\rm Law}}{=}  \{B_t,\,t\geq 0\}.
\end{equation}

For any choice of the Hurst parameter
$H\in(0,1)$, the Gaussian space generated by $B$ can be identified
with an isonormal Gaussian process of the type
$B=\{B(h):h\in\HH\}$, where the real and separable Hilbert space
$\EuFrak H$ is defined as follows: (i) denote by $\mathscr{E}$ the
set of all $\mathbb{R}$-valued step functions on $[0,\infty)$, (ii)
define $\EuFrak H$ as the Hilbert space obtained by closing
$\mathscr{E}$ with respect to the scalar product
$$
\left\langle
{\mathbf{1}}_{[0,t]},{\mathbf{1}}_{[0,s]}\right\rangle _{\EuFrak
H}=R(t,s).
$$
In particular, with such a notation, one has that
$B_t=B(\mathbf{1}_{[0,t]})$.

From now, assume on one hand that $B$ is defined on $[0,1]$ and on the other hand that 
the Hurst index verifies $H>\frac12$. The covariance kernel $R$ can be written as
\begin{equation*}
R(t,s)=\int_{0}^{s\wedge t}K_{H}(t,r)K_{H}(s,r)dr,
\end{equation*}%
where $K_{H}$ is the square integrable kernel given by
\begin{equation}\label{kh}
K_{H}(t,s)=\Gamma \left(H+\frac{1}{2}\right)^{-1}(t-s)^{H-\frac{1}{2}}F\left(H-\frac{1}{2},%
\frac{1}{2}-H,H+\frac{1}{2},1-\frac{t}{s}\right),
\end{equation}
$F(a,b,c,z)$ being the Gauss hypergeometric function. Consider the linear
operator $K_{H}^{\ast }$ from $\mathscr{E}$ to $L^{2}([0,1])$ defined by
\begin{equation*}
(K_{H}^{\ast }\varphi )(s)=K_{H}(1,s)\varphi (s)+\int_{s}^{1}\big( \varphi
(r)-\varphi (s)\big) \partial_1 K_{H}(r,s)dr.
\end{equation*}%
For any pair of step functions $\varphi $ and $\psi $ in $\mathscr{E}$, we
have
$\langle K_{H}^{\ast }\varphi ,K_{H}^{\ast }\psi \rangle
_{L^{2}}=\langle \varphi ,\psi \rangle _{\EuFrak{H}}$.
As a consequence, the operator $K_{H}^{\ast }$ provides an isometry between
the Hilbert spaces $\HH$ and $L^{2}([0,1])$. Hence, the process $%
W=(W_{t})_{t\in \lbrack 0,1]}$ defined by
\begin{equation}
W_{t}=B\left( (K_{H}^{\ast })^{-1}(\mathbf{1}_{[0,t]})\right)
\label{wiener}
\end{equation}%
is a Wiener process, and the process $B$ has an integral representation of
the form
$
B_{t}=\int_{0}^{t}K_{H}(t,s)dW_{s},
$
because $(K_{H}^{\ast }\mathbf{1}_{[0,t]})(s)=K_{H}(t,s)$.

The elements of $\HH$ may be not
functions but distributions. However, $\HH$ contains the subset $|\HH|$ of all measurable functions 
$f:[0,1]\to\mathbb{R}$ such that
$$
\int_{[0,1]^2}|f(u)||f(v)||u-v|^{2H-2}dudv<\infty.
$$
Moreover, for $f,g\in|\HH|$, we have
$$
\langle f,g\rangle_{\HH} = H(2H-1)\int_{[0,1]^2}f(u)\,g(v)\,|u-v|^{2H-2}dudv.
$$

In the sequel, we note $\HH^{\otimes q}$ and $\HH^{\odot q}$, respectively, the tensor space and the symmetric tensor space
of order $q\geq 1$. 
Let $\{e_k\::\: k\geq 1\}$ be a complete orthogonal system in $\HH$. 
Given $f\in\HH^{\odot p}$ and $g\in\HH^{\odot q}$, 
for every $r=0, \dots, p\wedge q$, the $r$th contraction of $f$ and $g$ is the element of 
$\HH^{\odot(p+q-2r)}$ defined as
$$
f\otimes_r g=\sum_{i_1=1, \dots, i_r=1}^{\infty} 
\langle f,e_{i_1}\otimes\cdots\otimes e_{i_r}\rangle_{\HH^{\otimes r}} 
\otimes
\langle g,e_{i_1}\otimes\cdots\otimes e_{i_r}\rangle_{\HH^{\otimes r}}. 
$$
In particular, note that $f\otimes_0g=f\otimes g$ and, when $p=q$, that $f\otimes_pg=\langle f, g\rangle_\HH$.
Since, in general, the contraction $f\otimes_rg$ is not a symmetric element of $\HH^{\otimes (p+q-2r)}$, we define 
$f\widetilde\otimes_rg$ as the canonical symmetrization of $f\otimes_rg$. 
When $f\in \HH^{\odot q}$, we write $I_q(f)$ to indicate its $q$th multiple integral with respect to $B$. 
The following formula is useful to compute the product of such integrals: 
if $f\in \HH^{\odot p}$ and $g\in \HH^{\odot q}$, then 
\begin{equation}
\label{eq:multiplication}
I_p(f)I_q(g)=\sum_{r=0}^{p\wedge q} r!
\left(\!\!\begin{array}{c}p\\r\end{array}\!\!\right)
\left(\!\!\begin{array}{c}q\\r\end{array}\!\!\right)
I_{p+q-2r}(f\widetilde\otimes_rg). 
\end{equation}

Let $\mathscr{S}$ be the set of cylindrical functionals $F$ of
the form
\begin{equation}
\label{eq:cylindrical} 
F=\varphi(B(h_1),\ldots,B(h_n)),
\end{equation}
where $n\geq 1$, $h_i \in \EuFrak{H}$ and the function $\varphi\in
\mathscr{C}^{\infty}(\real^n)$ is such that its partial
derivatives have polynomial growth. The Malliavin
derivative $DF$ of a functional $F$ of the form
(\ref{eq:cylindrical}) is the square integrable
$\EuFrak{H}$-valued random variable defined as
$$
DF=\sum_{i=1}^n \partial_i \varphi(B(h_1),\ldots,B(h_n)) h_i,
$$
where $\partial_i \varphi$ denotes the $i$th partial derivative of $\varphi$.
In particular, one has that $D_sB_t={\bf 1}_{[0,t]}(s)$ for every $s,t\in[0,1]$. 
As usual, $\mathbb{D}^{1,2}$ denotes the
closure of $\mathscr{S}$ with respect to the norm $\| \cdot
\|_{1,2}$, defined by the relation
$\| F\|_{1,2}^2 \; = \; E \big|F\big|^2 +
E \| D F\|_{\HH}^2.$
Note that every multiple integral belongs to $\mathbb{D}^{1,2}$.
Moreover, we have $$D_t\big(I_q(f)\big)=qI_{q-1}\big(f(\cdot,t)\big)$$ for any
$f\in\HH^{\odot q}$ and $t\geq 0$. 
The Malliavin derivative $D$ also satisfies the following
chain rule formula: if $\varphi:\real^n\to \real$ is continuously
differentiable with bounded derivatives and if
$(F_1,\ldots,F_n)$ is a random vector such that each component
belongs to $\mathbb{D}^{1,2}$, then $\varphi(F_1,\ldots,F_n)$ is
itself an element of $\mathbb{D}^{1,2}$, and moreover
$$
D\varphi(F_1,\ldots,F_n)=\sum_{i=1}^n \partial_i\varphi(F_1,\ldots,F_n) DF_i.
$$

\section{Case $H\in (1-1/(2q), 1)$}
\label{sec:H>}
In this section, we fix $q\geq 2$, we assume that $H>1-\frac1{2q}$
and we consider $Z$ defined by (\ref{her-rv}) for $W$ the Wiener
process defined by (\ref{wiener}).
By the scaling property (\ref{scaling}) of fBm, remark first 
that $Z_n$, defined by (\ref{eq:Breuer_Major3}), has the same law, for any fixed $n\geq 1$, as
\begin{equation}\label{Sn}
S_n=n^{q(1-H)-1}\sum_{k=0}^{n-1} H_q\big(n^{H}(B_{(k+1)/n}-B_{k/n})\big)=I_q(f_n),
\end{equation}
for $f_n=n^{q-1}\sum_{k=0}^{n-1} {\bf 1}_{[k/n,(k+1)/n]}^{\otimes q}\in\HH^{\odot q}$. 
In \cite{NNT}, Theorem 1 (point 3), it is shown that
the sequence $\{S_n\}_{n\geq 1}$ converges in $L^2$ towards $Z$, or equivalently that
$\{f_n\}_{n\geq 1}$ is Cauchy in $\HH^{\odot q}$.
Here, we precise the rate of this convergence: 
\begin{prop}
\label{thm-troisquart}
Let $f$ denote the limit of the Cauchy sequence $\{f_n\}_{n\geq 1}$ in $\HH^{\odot q}$. We have
$$
E\big|S_n - Z\big|^2 
=E\big|I_q(f_n) - I_q(f)\big|^2
= \|f_n - f\|^2_{\HH^{\odot q}} = O(n^{2q-1-2qH}),\quad\mbox{as $n\to\infty$}.
$$
\end{prop}
\noindent
Proposition \ref{thm-troisquart}, together with Theorem  \ref{theo:Davydov} above,
immediately entails \eqref{eq:H>} so that the rest of this section is devoted to the proof of 
the proposition.\\
\\
{\it Proof of Proposition \ref{thm-troisquart}}. 
We have
\begin{eqnarray}
\|f_n\|^2_{\HH^{\odot q}}&=&n^{2q-2}\,\sum_{k,l=0}^{n-1} \langle {\bf 1}_{[k/n,(k+1)/n]},{\bf 1}_{[l/n,(l+1)/n]}\rangle_\HH^q\notag\\
&=&H^q(2H-1)^q \,n^{2q-2}\,\sum_{k,l=0}^{n-1}\left(\int_{k/n}^{(k+1)/n} du \int_{l/n}^{(l+1)/n} dv |u-v|^{2H-2}\right)^q.\label{eq1}
\end{eqnarray}
By letting $n$ goes to infinity, we obtain
\begin{eqnarray}
\|f\|^2_{\HH^{\odot q}}&=&H^q(2H-1)^q\int_{[0,1]^2} |u-v|^{2qH-2q}dudv\notag\\
&=&H^q(2H-1)^q\,\sum_{k,l=0}^{n-1}\int_{k/n}^{(k+1)/n} du \int_{l/n}^{(l+1)/n} dv |u-v|^{2qH-2q}.\label{eq2}
\end{eqnarray}
Now, let $\phi\in|\HH|$. We have
\begin{eqnarray*}
\langle f_n,\phi^{\otimes q}\rangle_{\HH^{\odot q}}&=&\,n^{q-1}\,\sum_{l=0}^{n-1} \langle {\bf 1}_{[l/n,(l+1)/n]},\phi\rangle_\HH^q\\
&=&H^q(2H-1)^q \,n^{q-1}\,\sum_{l=0}^{n-1}\left(\int_{l/n}^{(l+1)/n} dv \int_{0}^{1} du \,\phi(u)|u-v|^{2H-2}\right)^q.
\end{eqnarray*}
By letting $n$ goes to infinity, we obtain
$$
\langle f,\phi^{\otimes q}\rangle_{\HH^{\odot q}} = 
H^q(2H-1)^q \int_0^1 dv \left(  \int_{0}^{1} du \,\phi(u)|u-v|^{2H-2}\right)^q.
$$
Hence, we have
\begin{eqnarray}
\langle f,f_n\rangle_{\HH^{\odot q}}&=&H^q(2H-1)^q \,n^{q-1}\,\sum_{k=0}^{n-1}\int_{0}^1 dv\left(\int_{k/n}^{(k+1)/n} du |u-v|^{2H-2}\right)^q\notag\\
&=&H^q(2H-1)^q \,n^{q-1}\,\sum_{k,l=0}^{n-1}\int_{l/n}^{(l+1)/n} dv\left(\int_{k/n}^{(k+1)/n} du |u-v|^{2H-2}\right)^q.\label{eq3}
\end{eqnarray}
Finally, by combining (\ref{eq1}), (\ref{eq2}) and (\ref{eq3}), and by using among others elementary change of variables, we can write:
\begin{eqnarray}
\nonumber
\|f_n-f\|^2_{\HH^{\odot q}}
&=&H^q(2H-1)^q \,\sum_{k,l=0}^{n-1}\left\{
n^{2q-2}\,\left(\int_{k/n}^{(k+1)/n} du \int_{l/n}^{(l+1)/n} dv |u-v|^{2H-2}\right)^q\right.\\
\nonumber
&&
-2\,n^{q-1}\,\int_{l/n}^{(l+1)/n} dv\left(\int_{k/n}^{(k+1)/n} du |u-v|^{2H-2}\right)^q\\
\nonumber
&&\left.
+\int_{k/n}^{(k+1)/n} dw \int_{l/n}^{(l+1)/n} dz |w-z|^{2qH-2q}
\right\}\\
\nonumber
&=&H^q(2H-1)^q \,n^{2q-2-2qH}\sum_{k,l=0}^{n-1}\left\{
\left(\int_{0}^{1} du \int_{0}^{1} dv |k-l+u-v|^{2H-2}\right)^q\right.\\
\nonumber
&&\left.
-2\,\int_{0}^{1} dv\left(\int_{0}^{1} du |k-l+u-v|^{2H-2}\right)^q
+\int_{0}^{1} du \int_{0}^{1} dv |k-l+u-v|^{2qH-2q}
\right\}\\
\nonumber
&\leq &H^q(2H-1)^q \,n^{2q-1-2qH}\sum_{r\in\mathbb{Z}}\left|
\left(\int_{0}^{1} du \int_{0}^{1} dv |r+u-v|^{2H-2}\right)^q\right.\\
\label{eq:sum_r}
&&\left.
-2\,\int_{0}^{1} dv\left(\int_{0}^{1} du |r+u-v|^{2H-2}\right)^q
+\int_{0}^{1} du \int_{0}^{1} dv |r+u-v|^{2qH-2q}
\right|.
\end{eqnarray}
Consequently, to achieve the proof of Theorem \ref{thm-troisquart}, it remains to ensure that 
the sum over $\mathbb{Z}$ in \eqref{eq:sum_r} is finite. For $r>1$, elementary computations give
\begin{eqnarray}
\nonumber
\left(\int_{0}^{1} du \int_{0}^{1} dv |r+u-v|^{2H-2}\right)^q
&=&\big(2H(2H-1)\big)^{-q}\big((r+1)^{2H}-2r^{2H}+(r-1)^{2H}\big)^q\\
\nonumber
&=&\left(r^{2H-2}+O(r^{2H-4})\right)^q\\
\label{eq:T1}
&=&r^{2qH-2q}+O(r^{2qH-2q-2})
\end{eqnarray}
and
\begin{eqnarray}
\int_{0}^{1} du \int_{0}^{1} dv |r+u-v|^{2qH-2q}
\nonumber
&=&\frac{(r+1)^{2qH-2q+2}-2r^{2qH-2q+2}+(r-1)^{2qH-2q+2}}{(2qH-2q+1)(2qH-2q+2)}\\
\label{eq:T3}
&=&r^{2qH-2q}+O(r^{2qH-2q-2}).
\end{eqnarray}
Moreover, using the inequality $\big|(1+x)^{2H-1}-1-(2H-1) x\big|\le (2H-1)(H-1) x^2$
for $x\in[0,1]$, we can write
\begin{eqnarray*}
\int_{0}^{1} dv
\left(\int_{0}^{1} du |r+u-v|^{2H-2}\right)^q
&=& (2H-1)^{-q}\int_{0}^{1} \big((r+1-v)^{2H-1}-(r-v)^{2H-1}\big)^qdv\\
&=& (2H-1)^{-q}\int_{0}^{1} (r-v)^{2qH-q}\left(\big(1+\frac{1}{r-v}\big)^{2H-1}-1\right)^qdv\\
&=& \int_{0}^{1} (r-v)^{2qH-q}\left(\frac 1{r-v}+R\big(\frac1{r-v}\big)\right)^qdv
\end{eqnarray*}
where the remainder term $R$ verifies
$
|R(u)|\leq(1-H)u^2.
$
In particular, for any $v\in[0,1]$, we have
$$
(r-v)\left|R\big(\frac 1{r-v}\big)\right|\leq \frac {1-H}{r-1}.
$$
Hence, we deduce:
\begin{eqnarray}
\nonumber
\int_{0}^{1} dv
\left(\int_{0}^{1} du |r+u-v|^{2H-2}\right)^q
&=& \int_{0}^{1} (r-v)^{2qH-2q}\left(1+O(1/r)\right)^qdv\\
\nonumber
&=& r^{2qH-2q+1}\frac{1-(1-1/r)^{2qH-2q+1}}{2qH-2q+1}(1+O(1/r))\\
\nonumber
&=& r^{2qH-2q+1}(1/r+O(1/r^2)(1+O(1/r))\\
\label{eq:T2}
&=& r^{2qH-2q}+O(r^{2qH-2q-1}).
\end{eqnarray}
By combining \eqref{eq:T1}, \eqref{eq:T3} and \eqref{eq:T2}, we obtain 
(since similar arguments also apply for $r<-1$) that
\begin{eqnarray*}
&&\left(\int_{0}^{1} du \int_{0}^{1} dv |r+u-v|^{2H-2}\right)^q
-2\,\int_{0}^{1} dv\left(\int_{0}^{1} du |r+u-v|^{2H-2}\right)^q\\
&&\hskip8cm+\int_{0}^{1} du \int_{0}^{1} dv |r+u-v|^{2qH-2q}
\end{eqnarray*}
is $O(|r|^{2qH-2q-1})$, 
so that the sum over $\mathbb{Z}$ in \eqref{eq:sum_r} is finite.  
The proof of Proposition \ref{thm-troisquart} is done.\qed

\section{Case $H=1-1/(2q)$}
\label{sec:H=}

As we already pointed out in the Introduction, the proof of \eqref{eq:H=} 
is a slight adaptation of that of Theorem \ref{theo:NP_Berry} (that is Theorem 4.1 in \cite{NP}) which was devoted to
the case when $H<1-1/(2q)$. 
That is why we will only focus, here, on the differences between the cases $H<1-1/(2q)$ and $H=1-1/(2q)$. 
In particular, we will freely refer to \cite{NP} each time we need an estimate already computed therein. 

From now, fix $H=1-1/(2q)$ and let us evaluate the right-hand side in Theorem \ref{theo:NP}.
Once again, instead of $Z_n$, we will rather use $S_n$ defined by 
$$
S_n=\frac{1}{\sigma_H\,\sqrt{n\log n}}\sum_{k=0}^{n-1} H_q\big(n^{H}(B_{(k+1)/n}-B_{k/n})\big)
=I_q\left(\frac{n^{q-1}}{\sigma_H\,\sqrt{\log n}}\sum_{k=0}^{n-1}
{\bf 1}_{[k/n,(k+1)/n]}^{\otimes q}
 \right)
$$ 
in the sequel, in order
to facilitate the connection with \cite{NP}.
First, observe that the covariance function $\rho_H$ of the Gaussian sequence 
$\big(n^H(B_{(r+1)/n}-B_{r/n})\big)_{r\geq 0}$,
given by 
$$
\rho_H(r)=\frac 12\big(|r+1|^{2-1/q}-2|r|^{2-1/q}+|r-1|^{2-1/q}\big),
$$ 
verifies the following straightforward expansion: 
\begin{equation}
\label{eq:dev_rho}
\rho_H(r)^q=\left((1-\frac 1{2q})(1-\frac 1q)\right)^q |r|^{-1}+O(|r|^{-3}),\quad\mbox{as $|r|\to\infty$}.
\end{equation}
Using $n^{2-1/q}\langle\ind_{[k/n,(k+1)/n]},\ind_{[l/n,(l+1)/n]}\rangle_{\HH}=\rho_H(k-l)$, note that
\begin{eqnarray*}
\Var(S_n)
&=&\frac{n^{2q-2}}{\sigma_H^2\log n}\sum_{k,l=0}^{n-1}\ee\big[I_q\big(\ind_{[k/n, (k+1)/n]}^{\otimes q}\big)I_q\big(\ind_{[l/n, (l+1)/n]}^{\otimes q}\big)\big]\\
&=&\frac{q!n^{2q-2}}{\sigma_H^2\log n}\sum_{k,l=0}^{n-1}\langle \ind_{[k/n, (k+1)/n]},\ind_{[l/n, (l+1)/n]}\rangle_\HH^q
=\frac{q!}{\sigma_H^2n\log n}\sum_{k,l=0}^{n-1}\rho_H(k-l)^q
\end{eqnarray*}
from which, together with \eqref{eq:dev_rho}, we deduce the exact value of $\sigma_H^2$: 
\begin{equation}
\label{eq:sigma}
\sigma_H^2:=\lim_{n\to +\infty}\frac{q!}{n\log n}\sum_{k,l=0}^{n-1}\rho_H(k-l)^q
=2q!\left(\frac{(2q-1)(q-1)}{2q^2}\right)^q.
\end{equation}

\noindent
In order to apply Theorem \ref{theo:NP}, we compute the Malliavin derivative of $S_n$: 
$$
DS_n=\frac{qn^{q-1}}{\sigma_H\sqrt{\log n}}\sum_{k=0}^{n-1}I_{q-1}\big(\ind_{[k/n, (k+1)/n]}^{\otimes q-1}\big)\ind_{[k/n,(k+1)/n]}.
$$
Hence
$$
\|DS_n\|_{\HH}^2\\
=\frac{q^2n^{2q-2}}{\sigma_H^2\log n}
\sum_{k,l=0}^{n-1}
I_{q-1}\big(\ind_{[k/n, (k+1)/n]}^{\otimes q-1}\big)I_{q-1}\big(\ind_{[l/n, (l+1)/n]}^{\otimes q-1}\big)
\langle\ind_{[k/n,(k+1)/n]},\ind_{[l/n,(l+1)/n]}\rangle_{\HH}.
$$
The multiplication formula \eqref{eq:multiplication} yields
\begin{align*}
&\|DS_n\|_{\HH}^2
=
\frac{q^2n^{2q-2}}{\sigma_H^2\log n}
\sum_{r=0}^{q-1}r! \left(\!\!\begin{array}{c} q-1\\r\end{array}\!\!\right)^2\times\\
&\hskip 2cm
\times \sum_{k,l=0}^{n-1}
I_{2q-2-2r}\big(\ind_{[k/n, (k+1)/n]}^{\otimes q-1-r} \tilde\otimes \ind_{[l/n, (l+1)/n]}^{\otimes q-1-r}\big)
\langle\ind_{[k/n,(k+1)/n]},\ind_{[l/n,(l+1)/n]}\rangle_{\HH}^{r+1}.
\end{align*}
We can rewrite 
$$
1-\frac 1q\|DS_n\|_{\HH}^2
=1-\sum_{r=0}^{q-1} A_r(n)
$$
where
\begin{eqnarray*}
A_r(n)&=&
\frac{q\,r!\left(\!\!\begin{array}{c} q-1\\r\end{array}\!\!\right)^2}{\sigma_H^2}\frac{n^{2q-2}}{\log n}\times \\
&&\times\sum_{k,l=0}^{n-1}I_{2q-2-2r}\big(\ind_{[k/n, (k+1)/n]}^{\otimes q-1-r} \tilde\otimes \ind_{[l/n, (l+1)/n]}^{\otimes q-1-r}\big)
\langle\ind_{[k/n,(k+1)/n]},\ind_{[l/n,(l+1)/n]}\rangle_{\HH}^{r+1}. 
\end{eqnarray*}
For the term $A_{q-1}(n)$, we have: 
\begin{eqnarray}
\nonumber
1-A_{q-1}(n)
&=&1-\frac {q!}{\sigma_H^2 \log n} n^{2q-2}\sum_{k,l=0}^{n-1} \langle\ind_{[k/n,(k+1)/n]},\ind_{[l/n,(l+1)/n]}\rangle_{\HH}^q\\
\nonumber
&=&1-\frac {q!}{\sigma_H^2 n\log n}\sum_{k,l=0}^{n-1} \rho_H(k-l)^q
=1-\frac {q!}{\sigma_H^2 n\log n}\sum_{|r|<n}(n-|r|) \rho_H(r)^q\\
\nonumber
&=&1-\frac {q!}{\sigma_H^2\log n}\sum_{|r|<n} \rho_H(r)^q+
\frac {q!}{\sigma_H^2 n\log n}\sum_{|r|<n}|r| \rho_H(r)^q
\label{eq:bound1}
=O(1/\log n)
\end{eqnarray}
where the last estimate comes from the development \eqref{eq:dev_rho} of $\rho_H$ 
and from the exact value \eqref{eq:sigma} of $\sigma_H^2$.

\medskip
\noindent
Next, we show that for any {\sl fixed} $r\leq q-2$, we have 
$\ee|A_r(n)|^2=O(1/\log n)$. Indeed: 
\begin{align*}
&\ee|A_r(n)|^2\\
&= c(H,r,q)\frac{n^{4q-4}}{\log^2 n} \sum_{i,j,k,l=0}^{n-1}
\langle \ind_{[k/n,(k+1)/n]}, \ind_{[l/n, (l+1)/n]}\rangle_{\HH}^{r+1}
\langle \ind_{[k/n,(k+1)/n]}, \ind_{[l/n, (l+1)/n]}\rangle_{\HH}^{r+1}\\
&\times\langle \ind_{[k/n,(k+1)/n]}^{\otimes q-1-r}\tilde\otimes\ind_{[l/n, (l+1)/n]}^{\otimes q-1-r}, 
 \ind_{[i/n,(i+1)/n]}^{\otimes q-1-r}\otimes\ind_{[j/n, (j+1)/n]}^{\otimes q-1-r}\rangle_{\HH^{\otimes 2q-2-2r}}\\
&=\sum_{\begin{subarray}{c}\gamma, \delta\geq 0\\\gamma+\delta=q-r-1\end{subarray}}
\sum_{\begin{subarray}{c}\alpha, \beta\geq 0\\\alpha+\beta=q-r-1\end{subarray}}
c(H,r,q,\alpha,\beta,\gamma,\delta) B_{r,\alpha,\beta,\gamma,\delta}(n)
\end{align*}
where $c(\cdot)$ is a generic constant depending only on its arguments and 
\begin{align*}
&B_{r,\alpha,\beta,\gamma,\delta}(n)
=\frac{n^{4q-4}}{\log n^2}\sum_{i,j,k,l=0}^{n-1}\\
&
\langle \ind_{[k/n,(k+1)/n]}, \ind_{[l/n, (l+1)/n]}\rangle_{\HH}^{r+1}
\langle \ind_{[i/n,(i+1)/n]}, \ind_{[j/n, (j+1)/n]}\rangle_{\HH}^{r+1}
\langle \ind_{[k/n,(k+1)/n]}, \ind_{[i/n, (i+1)/n]}\rangle_{\HH}^{\alpha}\\
&\langle \ind_{[k/n,(k+1)/n]}, \ind_{[j/n, (j+1)/n]}\rangle_{\HH}^{\beta}
\langle \ind_{[l/n,(l+1)/n]}, \ind_{[i/n, (i+1)/n]}\rangle_{\HH}^{\gamma}
\langle \ind_{[l/n,(l+1)/n]}, \ind_{[j/n, (j+1)/n]}\rangle_{\HH}^{\delta}\\
&=\frac{1}{n^2\log^2 n}\sum_{i,j,k,l=0}^{n-1}\rho_H(k-l)^{r+1}\rho_H(i-j)^{r+1}
\rho_H(k-i)^{\alpha}\rho_H(k-j)^{\beta}\rho_H(l-i)^{\gamma}\rho_H(l-j)^{\delta}.
\end{align*}
As in \cite{NP}, when $\alpha,\beta,\gamma,\delta$ are fixed, we can decompose the sum appearing in
$B_{r,\alpha,\beta,\gamma,\delta}(n)$ as follows: 
\begin{align*}
&\sum_{i=j=k=l}+
\left(
\sum_{\begin{subarray}{c}i=j=k\\l\not =i\end{subarray}}
+\sum_{\begin{subarray}{c}i=j=l\\k\not =i\end{subarray}}
+\sum_{\begin{subarray}{c}i=k=l\\j\not =i\end{subarray}}
+\sum_{\begin{subarray}{c}j=k=l\\i\not =j\end{subarray}}
\right)
+\left(
\sum_{\begin{subarray}{c}i=j, k=l\\k\not =i\end{subarray}}
+\sum_{\begin{subarray}{c}i=k, j=l\\j\not =i\end{subarray}}
+\sum_{\begin{subarray}{c}i=l, j=k\\j\not =i\end{subarray}}
\right)
\\
&+\left(
\sum_{\begin{subarray}{c}i=j, k\not=i\\k\not =l, l\not=i\end{subarray}}
+\sum_{\begin{subarray}{c}i=k, j\not=i\\j\not =l, k\not=l\end{subarray}}
+\sum_{\begin{subarray}{c}i=l, k\not=i\\k\not =j, j\not=i\end{subarray}}
+\sum_{\begin{subarray}{c}j=k, k\not=i\\k\not =l, l\not=i\end{subarray}}
+\sum_{\begin{subarray}{c}j=l, k\not=i\\k\not =l, l\not=i\end{subarray}}
+\sum_{\begin{subarray}{c}k=l, k\not=i\\k\not =j, j\not=i\end{subarray}}
\right)
+\sum_{\begin{subarray}{c}i,j,k,l\\ \mbox{{\footnotesize are all distinct}}\end{subarray}}
\end{align*}
where the indices run over $\{0, 1\dots, n-1\}$. 
Similar computations as in \cite{NP} show that 
the first, second and third sums are $O(1/(n\log^2 n))$; 
the fourth and  fifth sums are $O(1/(n^{1/q}\log^2 n))$; 
the sixth, seventh and eihght sums are $O(1/(n^{2/q}\log^2 n))$;
the ninth, tenth, eleventh, twelfth, thirteenth and fourteenth sums are also  $O(1/(n^{1/q}\log^2 n))$.
We focus only on the last sum. Actually, 
it is precisely its  
contribution which will indicate 
the true order in \eqref{eq:H=}.
Once again, we split this sum into $24$ sums:
\begin{equation}
\label{eq:15}
\sum_{k>l>i>j}+\sum_{k>l>j>i}+\cdots
\end{equation}
We first deal with the first one, for which we have 
\begin{align*}
&\frac 1{(n\log n)^2}\sum_{k>l>i>j}
\rho_H(k-l)^{r+1}\rho_H(i-j)^{r+1}\rho_H(k-i)^{\alpha}
\rho_H(k-j)^{\beta}\rho_H(l-i)^{\gamma}\rho_H(l-j)^{\delta}\\
&\trianglelefteqslant\frac 1{(n\log n)^2}\sum_{k>l>i>j}
(k-l)^{-1}(i-j)^{-(r+1)/q}(l-i)^{-(q-r-1)/q}\\
&=\frac 1{(n\log n)^2}\sum_k\sum_{l<k}(k-l)^{-1}\sum_{i<l} (l-i)^{-(q-r-1)/q}
\sum_{j<i}(i-j)^{-(r+1)/q}\\
&\trianglelefteqslant\frac 1{n\log^2 n}
\sum_{l=1}^{n-1}l^{-1}
\sum_{i=1}^{n-1} i^{-(q-r-1)/q}
\sum_{j=1}^{n-1}j^{-(r+1)/q}
\trianglelefteqslant\frac 1{n\log^2 n}
\sum_{l=1}^{n-1}l^{-1}
n^{(r+1)/q}n^{1-(r+1)/q}
\trianglelefteqslant \frac 1{\log n}
\end{align*}
where the notation $a_n\trianglelefteqslant b_n$ means that $\sup_{n\geq 1}|a_n|/|b_n|<+\infty$. 
Since the other sums in \eqref{eq:15} are similarly bounded, the fifteenth sum is $O(1/\log n)$. 
Consequently: 
$$
\sum_{r=0}^{q-2} \ee[A_r(n)^2]=O(1/\log n). 
$$
Finally, together with \eqref{eq:bound1}, we obtain
$
\ee\big[\big(1-\frac 1q\|DS_n\|_\HH^2\big)^2\big]
= O(1/\log n)
$
and the proof of \eqref{eq:H=} is achieved thanks to Theorem \ref{theo:NP}.

\end{document}